\documentclass[11pt,a4paper]{article}

\usepackage[T1]{fontenc}
\usepackage[utf8]{inputenc}
\usepackage{lmodern}
\usepackage{amsmath,amssymb,amsthm,amsfonts}
\usepackage{mathrsfs}
\usepackage{geometry}
\usepackage{enumitem}
\usepackage{cite}
\usepackage[hidelinks]{hyperref}

\theoremstyle{plain}
\newtheorem{theorem}{Theorem}[section]
\newtheorem{lemma}[theorem]{Lemma}
\newtheorem{proposition}[theorem]{Proposition}
\newtheorem{corollary}[theorem]{Corollary}

\theoremstyle{definition}
\newtheorem{definition}[theorem]{Definition}
\newtheorem{example}[theorem]{Example}

\theoremstyle{remark}
\newtheorem{remark}[theorem]{Remark}

\numberwithin{equation}{section}

\title{On Partially Totally Real Riemannian Maps and Their Sharp Curvature Inequalities}
\author{Tanveer Fatima}
\date{}

\begin{document}

\maketitle

\begin{abstract}
We introduce partially totally real Riemannian maps into K\"ahler manifolds and examine the geometry of their associated distributions. For complex space forms, we obtain a Chen-type curvature inequality retaining the contributions of the mixed second fundamental form and O'Neill's tensor, together with an arbitrary-codimension estimate involving partial mean curvatures. We characterize equality in the principal inequality and give an example attaining the bound when $\dim\mathcal D\geq4$, including a case with nonintegrable horizontal distribution.
\end{abstract}

\begingroup
\renewcommand{\thefootnote}{}
\footnotetext{%


\medskip
\noindent\textbf{Author's address:}
Department of Mathematics and Statistics,
College of Sciences in Yanbu, Taibah University,
Yanbu Governorate, Kingdom of Saudi Arabia.

\noindent\textbf{E-mail address:}
\texttt{tansari@taibahu.edu.sa}}
\endgroup

\section{Introduction}

Riemannian maps, introduced by Fischer~\cite{Fischer1992}, provide a common generalization of isometric immersions and Riemannian submersions. Their geometry is determined by the interaction between the horizontal and vertical distributions of the domain and the range and normal distributions along the map. The fundamental equations of Riemannian maps have been developed in this setting, leading to the study of various geometric properties, including harmonicity and totally geodesic maps.

When the target is an almost Hermitian manifold, the action of the almost complex structure on the range distribution provides further geometric structure. In this direction, \c{S}ahin~\cite{Sahin2010,Sahin2011,Sahin2017Generic} investigated invariant, anti-invariant, semi-invariant and generic Riemannian maps. In particular, semi-invariant Riemannian maps admit an orthogonal decomposition of the range distribution into invariant and anti-invariant components. Related range-side slant, semi-slant and hemi-slant structures have also been studied in~\cite{Sahin2013Slant,ParkSahin2014,Sahin2017HemiSlant,GunduzalpAkyolSahin2024}. These investigations motivate the consideration of more general decompositions without imposing invariance or slant conditions on the complementary distribution.

A decomposition of this kind underlies the partially totally real (PTR) submanifolds introduced by Poyraz, \c{S}ahin and Yerlikaya~\cite{PoyrazSahinYerlikaya2025}. Its geometric content is not restricted to the behaviour of $J$: the integrability and total geodesicity of the two distributions depend on the normal second fundamental form. For Riemannian maps, one must also account for the vertical components of covariant derivatives. Thus the corresponding integrability and local product conditions naturally involve both the shape operators of the range and O'Neill's tensor $\mathcal A$.

Curvature estimates present a related distinction. Given the above decomposition, the difference between the horizontal scalar-curvature trace and its restriction to $\mathcal D$ records precisely those sectional curvatures whose two-planes meet $\mathcal D^\perp$. In a complex space form, the ambient sectional curvature of every such plane is $c/4$. The Gauss--O'Neill equation then separates the ambient contribution from terms determined by the second fundamental form and the nonintegrability of the horizontal distribution. For PTR submanifolds, a corresponding distribution-dependent invariant and curvature inequality were considered by Mihai et al.~\cite{MihaiEtAl2026}.
Among the existing results are Chen's inequality for
Riemannian maps obtained by \c{S}ahin~\cite{Sahin2016Chen}
and the corresponding curvature estimates for maps into
complex space forms established by Meena et al.~\cite{MeenaSahinShah2025}. Further Chen-type and Casorati inequalities for various
classes of Riemannian maps and Riemannian submersions
have been studied in
~\cite{FatimaEtAl2025,SinghMeenaMeena2026,
Singh2026Generalized,FatimaAlDayel2026,
FatimaAlharbiDeshmukh2026}.
.

In this paper we define PTR Riemannian maps without imposing
maximality on $\mathcal D^\perp$ and examine the geometry of
the resulting distributions. Section~2 recalls the basic notions and fundamental equations needed in the sequel.
In Section~3, we establish integrability and totally geodesic
conditions for the associated horizontal distributions and
characterize their local product structure. Section~4 is
devoted to curvature inequalities for PTR Riemannian maps
into complex space forms. We consider the horizontal
curvature difference associated with the chosen PTR
decomposition. The K\"ahler structure yields a symmetric
cubic form on the totally real part of the range; combining
its algebraic estimate with the horizontal Gauss--O'Neill
equation gives a Chen-type inequality that retains the
contributions of the mixed second fundamental form and
O'Neill's tensor $\mathcal A$. We determine the equality
conditions and construct an example attaining the bound
with nonzero mean curvature and nonintegrable horizontal
distribution. An additional estimate in arbitrary normal
codimension retains the partial mean-curvature contributions
and makes the relation with the PTR-submanifold inequality
explicit under minimum codimension. Equality in the
principal bound does not require mixed total geodesicity.

\section{Preliminaries}

Let $(N,g_N,J)$ be an almost Hermitian manifold
with almost complex structure $J$ satisfying
\begin{equation}\label{eq:almost_hermitian}
J^2=-I,\qquad g_N(JU,JV)=g_N(U,V)
\end{equation}
for all $U,V\in\Gamma(TN)$. The fundamental
$2$-form $\Omega$ is defined by
$\Omega(U,V)=g_N(JU,V)$. An almost Hermitian
manifold is called K\"ahler if
\begin{equation}\label{eq:kahler}
(\nabla^N_U J)V=0,
\end{equation}
where $\nabla^N$ denotes the Levi-Civita
connection on $N$.

A K\"ahler manifold of constant holomorphic sectional curvature $c$ is called a complex space form and is denoted by $N(c)$. With the convention $K_N(U,V)=g_N(R^N(U,V)V,U)$ for orthonormal $U,V$, its curvature tensor is given by
\begin{align}\label{eq:complex_space_form}
R^N(U,V)W=\frac{c}{4}\big\{
&g_N(V,W)U-g_N(U,W)V \notag\\
&+g_N(JV,W)JU-g_N(JU,W)JV \notag\\
&+2g_N(U,JV)JW
\big\}.
\end{align}

Let $F:(M,g_M)\to(N,g_N)$ be a smooth map
of constant rank. The map $F$ is called
a Riemannian map if
\begin{equation}\label{eq:riemannian_map}
g_N(F_*X,F_*Y)=g_M(X,Y)
\end{equation}
for all
$X,Y\in\Gamma((\ker F_*)^\perp)$.

We denote the vertical and horizontal
distributions by
\[
\mathcal V=\ker F_*, \qquad \mathcal H=(\ker F_*)^\perp,
\]
respectively. Thus,
\[
TM=\mathcal V\oplus\mathcal H.
\]
Along $F$, we also have
\[
F^{-1}TN =\operatorname{range}F_* \oplus(\operatorname{range}F_*)^\perp.
\]

Let $\nabla^F$ be the pullback connection
induced by $\nabla^N$. The second
fundamental form of $F$ is defined by~\cite{Nore1987}
\begin{equation}\label{eq:second_fundamental_form}
(\nabla F_*)(E,G) =\nabla^F_EF_*G-F_*(\nabla^M_EG),
\end{equation}
where $\nabla^M$ is the Levi-Civita
connection on $M$.

For horizontal vector fields $X,Y$,
the second fundamental form satisfies
\begin{equation}\label{eq:normal_second_form}
(\nabla F_*)(X,Y) \in\Gamma((\operatorname{range}F_*)^\perp).
\end{equation}
For this standard property of Riemannian maps, see~\cite{Sahin2017Book}. Consequently, the Gauss formula for
a Riemannian map is
\begin{equation}\label{eq:gauss_map}
\nabla^F_XF_*Y =F_*(\nabla^M_XY)+(\nabla F_*)(X,Y).
\end{equation}

For a normal vector field
$\xi\in\Gamma((\operatorname{range}F_*)^\perp)$,
the Weingarten formula is

the Weingarten formula is
\begin{equation}\label{eq:weingarten_map}
\nabla^F_X\xi =-\mathcal S_\xi F_*X +\nabla_X^{F,\perp}\xi,
\end{equation}
where $\nabla^{F,\perp}$ is the normal connection. We write $\mathcal S_\xi$ for the shape operator, reserving $\mathcal A$ for O'Neill's tensor; thus
\begin{equation}\label{eq:shape_operator}
g_N(\mathcal S_\xi F_*X,F_*Y) =g_N(\xi,(\nabla F_*)(X,Y)).
\end{equation}

We recall the fundamental tensors
$\mathcal T$ and $\mathcal A$ associated
with the vertical and horizontal
distributions:
\begin{align}
\mathcal T_EG
&=\mathcal H\nabla^M_{\mathcal VE}\mathcal VG
+\mathcal V\nabla^M_{\mathcal VE}\mathcal HG,
\label{eq:t_tensor}\\
\mathcal A_EG
&=\mathcal H\nabla^M_{\mathcal HE}\mathcal VG
+\mathcal V\nabla^M_{\mathcal HE}\mathcal HG.
\label{eq:a_tensor}
\end{align}
Here $\mathcal V$ and $\mathcal H$ also
denote the corresponding orthogonal
projections.

For $U,V\in\Gamma(\mathcal V)$ and
$X,Y\in\Gamma(\mathcal H)$, these tensors
give
\begin{align}
\nabla^M_UV
&=\mathcal T_UV+\mathcal V\nabla^M_UV,
\label{eq:vertical_connection}\\
\nabla^M_UX
&=\mathcal H\nabla^M_UX+\mathcal T_UX,
\label{eq:mixed_connection_one}\\
\nabla^M_XU
&=\mathcal A_XU+\mathcal V\nabla^M_XU,
\label{eq:mixed_connection_two}\\
\nabla^M_XY
&=\mathcal H\nabla^M_XY+\mathcal A_XY.
\label{eq:horizontal_connection}
\end{align}

In particular,
\begin{equation}\label{eq:a_bracket}
\mathcal A_XY=\frac12\mathcal V[X,Y].
\end{equation}
This identity follows from the fundamental equations of a Riemannian submersion~\cite{ONeill1966,Sahin2017Book}. Thus, the horizontal distribution is
integrable if and only if
$\mathcal A_XY=0$ for all horizontal
vector fields $X,Y$.

\section{PTR Riemannian Maps}

In this section, we introduce partially totally
real Riemannian maps into K\"ahler manifolds
and establish some basic properties of the
distributions associated with these maps.

\begin{definition}\label{def:ptr_map}
Let $F:(M,g_M)\to(N,g_N,J)$ be a Riemannian
map into a K\"ahler manifold. The map $F$
is called a \emph{partially totally real
Riemannian map}, or briefly a PTR Riemannian
map, if $\operatorname{range}F_*$ admits an orthogonal decomposition
\begin{equation}\label{eq:ptr_decomposition}
\operatorname{range}F_*=\mathcal D\oplus\mathcal D^\perp,
\end{equation}
where $\mathcal D^\perp$ is a nonzero totally real distribution of constant rank, that is,
\begin{equation}\label{eq:ptr_totally_real}
J\mathcal D^\perp \subseteq(\operatorname{range}F_*)^\perp,
\end{equation}
and the tangential component of $J$
preserves $\mathcal D$.
\end{definition}

Following the terminology for PTR submanifolds, $\mathcal D$ will be called the ambiguous (or complementary) distribution. We call $F$ a \emph{proper PTR Riemannian map} if both distributions in \eqref{eq:ptr_decomposition} are nonzero and $\mathcal D$ is neither invariant, totally real, slant nor pointwise slant. 

In the above definition, no maximality condition
is imposed on $\mathcal D^\perp$; thus, it need
not coincide with $\ker P$.

For any $F_*X\in\operatorname{range}F_*$,
we write
\begin{equation}\label{eq:ptr_J_tangent}
JF_*X=PF_*X+QF_*X,
\end{equation}
where $PF_*X\in\operatorname{range}F_*$
and $QF_*X\in(\operatorname{range}F_*)^\perp$.

Similarly, for
$\xi\in(\operatorname{range}F_*)^\perp$,
we put
\begin{equation}\label{eq:ptr_J_normal}
J\xi=t\xi+f\xi,
\end{equation}
where $t\xi\in\operatorname{range}F_*$
and $f\xi\in(\operatorname{range}F_*)^\perp$.

Since $P$ is skew-adjoint, the totally real condition gives
\begin{equation}\label{eq:ptr_P_properties}
P\mathcal D^\perp=0, \qquad P\mathcal D\subseteq\mathcal D.
\end{equation}
Indeed, for $U\in\mathcal D$ and $W\in\mathcal D^\perp$, we have
\[
g_N(PU,W) =-g_N(U,JW)=0.
\]

\begin{lemma}\label{lem:ptr_normal_decomposition}
Let $F:(M,g_M)\to(N,g_N,J)$ be a PTR
Riemannian map. Suppose that
$Q|_{\mathcal D}$ has constant rank.
Then the normal bundle admits the
orthogonal decomposition
\begin{equation}\label{eq:ptr_normal_decomposition}
(\operatorname{range}F_*)^\perp =J\mathcal D^\perp \oplus Q\mathcal D \oplus\mu,
\end{equation}
where $\mu$ is invariant under $J$.
\end{lemma}

\begin{proof}
For $U\in\mathcal D$ and
$W\in\mathcal D^\perp$, we obtain
\[
g_N(JW,QU) =g_N(JW,JU) =g_N(W,U)=0.
\]
Thus $J\mathcal D^\perp$ and
$Q\mathcal D$ are orthogonal.

Let $\mu$ be their orthogonal complement
in $(\operatorname{range}F_*)^\perp$.
For $\xi\in\mu$, we have
\[
g_N(J\xi,U) =-g_N(\xi,QU)=0,
\]
and
\[
g_N(J\xi,W) =-g_N(\xi,JW)=0.
\]
Hence $J\xi$ is normal to
$\operatorname{range}F_*$.

Moreover,
\[
g_N(J\xi,JW)=g_N(\xi,W)=0
\]
and
\[
g_N(J\xi,QU) =g_N(J\xi,JU-PU)=0.
\]
Consequently, $J\xi\in\mu$.
Since $J^2=-I$, it follows that
$J\mu=\mu$.
\end{proof}

The following relations are direct
consequences of the preceding
decompositions and $J^2=-I$.

For $W\in\mathcal D^\perp$,
\begin{equation}\label{eq:ptr_algebraic_one}
PW=0,\qquad QW=JW,\qquad t(JW)=-W,\qquad f(JW)=0.
\end{equation}

For $U\in\mathcal D$,
\begin{equation}\label{eq:ptr_algebraic_two}
t(QU)=-U-P^2U,\qquad f(QU)=-Q(PU).
\end{equation}

Finally, for $\xi\in\mu$,
\begin{equation}\label{eq:ptr_algebraic_three}
t\xi=0,\qquad f\xi=J\xi.
\end{equation}

\begin{remark}\label{rem:ptr_special_cases}
The PTR structure encompasses several familiar classes of Riemannian maps, provided that the totally real distribution is nonzero.

\begin{enumerate}[label=(\roman*)]
\item Every semi-invariant Riemannian map with a nonzero
anti-invariant distribution is PTR, with
$J\mathcal D=\mathcal D$. Anti-invariant Riemannian maps
correspond to $\mathcal D=\{0\}$, whereas invariant
Riemannian maps correspond formally to
$\mathcal D^\perp=\{0\}$, which is excluded from
our definition.

\item Hemi-slant and pointwise hemi-slant Riemannian
maps into almost Hermitian manifolds are PTR whenever
their range distributions have nonzero totally real
components. In these cases, $\mathcal D$ is the slant
or pointwise slant distribution, respectively
~\cite{GunduzalpAkyolSahin2024}.

\item A generic Riemannian map in the sense of
\cite{Sahin2017Generic} need not be PTR. Indeed,
allowing $\mathcal D^\perp=\{0\}$ would make the
PTR condition vacuous. However, a PTR decomposition
exists whenever $\ker P$ contains a nonzero smooth
subbundle of constant rank.

\item The isometric immersion of a PTR-submanifold
is a PTR Riemannian map with $\ker F_*=\{0\}$.
In particular, totally real, proper CR-, and
hemi-slant submanifolds with nonzero totally real
components provide examples.

\item Let $\pi:(M,g_M)\to(L,g_L)$ be a Riemannian
submersion and $i:(L,g_L)\to(N,g_N,J)$ an isometric
immersion of a PTR-submanifold. Then
$F=i\circ\pi$ is a PTR Riemannian map.
\end{enumerate}
The PTR decomposition is part of the specified geometric structure; without a maximality condition it need not be unique. In particular, the results below are relative to the chosen distributions.
\end{remark}

We now give two examples of proper PTR
Riemannian maps which are not included
in the preceding special classes.

\begin{example}\label{ex:ptr_nonintegrable}
Let $\mathbb C^5$ be equipped with its
standard K\"ahler structure. Choose
$0<\theta_1,\theta_2<\pi/2$ such that
$\theta_1\ne\theta_2$, and define
\[
\begin{aligned}
U_1&=\partial_{x_1},\\
U_2&=\cos\theta_1\,\partial_{y_1}
     +\sin\theta_1\,\partial_{x_2},\\
U_3&=\partial_{x_3},\\
U_4&=\cos\theta_2\,\partial_{y_3}
     +\sin\theta_2\,\partial_{x_4},\\
U_5&=\partial_{x_5}.
\end{aligned}
\]
These vectors form an orthonormal set.

Consider $\mathbb R^6$ with coordinates
$(u_1,\ldots,u_5,v)$ and the metric
\[
g_M=\sum_{i=1}^5du_i^2 +(dv+u_1\,du_5)^2.
\]
Define
\[
F:(\mathbb R^6,g_M)\longrightarrow \mathbb C^5,\qquad F(u_1,\ldots,u_5,v)=\sum_{i=1}^5u_iU_i.
\]

Then
\[
\ker F_*=\operatorname{span}\{\partial_v\},
\]
and an orthonormal frame of
$(\ker F_*)^\perp$ is
\[
X_i=\partial_{u_i},\quad 1\leq i\leq4, \qquad X_5=\partial_{u_5}-u_1\partial_v.
\]
Since $F_*X_i=U_i$, the map $F$ is
Riemannian.

Set
\[
\mathcal D=\operatorname{span} \{U_1,U_2,U_3,U_4\}, \qquad \mathcal D^\perp=\operatorname{span}\{U_5\}.
\]
Clearly,
\[
J\mathcal D^\perp \subseteq(\operatorname{range}F_*)^\perp.
\]
Moreover,
\[
\begin{aligned}
PU_1&=\cos\theta_1U_2,
&PU_2&=-\cos\theta_1U_1,\\
PU_3&=\cos\theta_2U_4,
&PU_4&=-\cos\theta_2U_3.
\end{aligned}
\]
Thus
\[
P^2|_{\operatorname{span}\{U_1,U_2\}} =-\cos^2\theta_1 I,
\]
whereas
\[
P^2|_{\operatorname{span}\{U_3,U_4\}} =-\cos^2\theta_2 I.
\]
Since $\theta_1\ne\theta_2$,
$\mathcal D$ is neither invariant
nor pointwise slant. Hence $F$ is
a proper PTR Riemannian map.

Finally,
\[
[X_1,X_5]=-\partial_v,
\]
so that
\[
\mathcal A_{X_1}X_5 =-\frac12\partial_v\ne0.
\]
Therefore, the horizontal distribution
of $F$ is nonintegrable.
\end{example}

\begin{example}\label{ex:ptr_nonzero_second_form}
Choose $0<a_1<a_2<1$ and put $b_j=\sqrt{1-a_j^2}$, $j=1,2$.
Choose a nonzero constant $k$ and a
sufficiently small interval $I$
containing $0$.

Consider
\[
M=I\times\mathbb R^5
\]
with coordinates $(s,u,p,z,w,v)$
and metric
\[
\begin{aligned}
g_M={}&(1+4k^2s^2)\,ds^2+du^2
       +dp^2+dz^2+dw^2\\
     &+(dv+s\,dw)^2.
\end{aligned}
\]

Define $F:M\to\mathbb C^5$ by
\[
F(s,u,p,z,w,v) = \bigl(s+ia_1u,\, b_1u+iks^2,\, p+ia_2z,\, b_2z,\, w\bigr).
\]

Put
\[
\rho(s)=\sqrt{1+4k^2s^2}.
\]
The vertical distribution is
\[
\ker F_*=\operatorname{span}\{\partial_v\},
\]
and the horizontal distribution has
the orthonormal frame
\[
\begin{aligned}
X_1&=\rho(s)^{-1}\partial_s,\\
X_2&=\partial_u,\qquad
X_3=\partial_p,\qquad
X_4=\partial_z,\\
X_5&=\partial_w-s\partial_v.
\end{aligned}
\]
Direct computation gives
\[
g_N(F_*X_i,F_*X_j)=\delta_{ij}.
\]
Hence $F$ is a Riemannian map.

Let
\[
\mathcal D=\operatorname{span} \{F_*X_1,F_*X_2,F_*X_3,F_*X_4\},
\]
and
\[
\mathcal D^\perp =\operatorname{span}\{F_*X_5\}.
\]
Since $F_*X_5=\partial_{x_5}$,
we have
\[
JF_*X_5=\partial_{y_5} \in(\operatorname{range}F_*)^\perp.
\]

The tangential component $P$ preserves
the two orthogonal distributions
\[
\mathcal D_1=\operatorname{span} \{F_*X_1,F_*X_2\}, \qquad \mathcal D_2=\operatorname{span} \{F_*X_3,F_*X_4\}.
\]
On these distributions,
\[
P^2|_{\mathcal D_1} =-\frac{(a_1-2b_1ks)^2}{1+4k^2s^2}\,I,
\]
and
\[
P^2|_{\mathcal D_2}= -a_2^2 I.
\]
By restricting $I$ if necessary,
we may assume that
\[
0< \frac{a_1-2b_1ks}{\rho(s)}<a_2<1.
\]
Consequently, $\mathcal D$ is
neither invariant nor pointwise slant.
Thus $F$ is a proper PTR Riemannian map.

Moreover,
\[
[X_1,X_5] =-\rho(s)^{-1}\partial_v,
\]
which implies
\[
\mathcal A_{X_1}X_5 =-\frac{1}{2\rho(s)}\partial_v \ne0.
\]

The second fundamental form of $F$
is also nonzero. Indeed, at $s=0$,
\[
(\nabla F_*)(X_1,X_1) =2k\,\partial_{y_2}\ne0.
\]
Thus this example has both a
nonintegrable horizontal distribution
and a nonzero horizontal second
fundamental form.
\end{example}

\subsection{Geometry of the PTR Distributions}

We study the horizontal distributions whose images under $F_*$ are $\mathcal D$ and $\mathcal D^\perp$, respectively. Unless stated otherwise, all vector fields in the following assertions are taken locally.

The following symmetry, corresponding to the shape-operator relation in the submanifold setting~\cite{PoyrazSahinYerlikaya2025}, will be used in the curvature estimate.

\begin{lemma}\label{lem:ptr_cubic}
Let $F:(M,g_M)\to(N,g_N,J)$ be a PTR Riemannian map. If $X,W,Z\in\Gamma(\mathcal H)$ and $F_*W,F_*Z\in\mathcal D^\perp$, then
\begin{equation}\label{eq:ptr_cubic}
g_N((\nabla F_*)(X,W),JF_*Z)=g_N((\nabla F_*)(X,Z),JF_*W).
\end{equation}
Consequently, the trilinear expression $g_N((\nabla F_*)(W,Z),JF_*U)$ is symmetric in $W,Z,U$ whenever their images lie in $\mathcal D^\perp$.
\end{lemma}

\begin{proof}
Since $\mathcal D^\perp$ is totally real, $g_N(F_*W,JF_*Z)=0$. Differentiating along $X$, using $\nabla^NJ=0$ and the Gauss formula for Riemannian maps, we obtain
\[
0=g_N((\nabla F_*)(X,W),JF_*Z)-g_N((\nabla F_*)(X,Z),JF_*W).
\]
Indeed, the terms containing $F_*(\nabla_X^MW)$ and $F_*(\nabla_X^MZ)$ vanish because $JF_*W$ and $JF_*Z$ are normal to $\operatorname{range}F_*$. This gives \eqref{eq:ptr_cubic}. The last assertion follows by combining this identity with the symmetry of $(\nabla F_*)$ on horizontal vectors.
\end{proof}

For the next assertion, assume that the totally real summand is maximal in the sense that $\ker P=\mathcal D^\perp$. Then $F_*^{-1}(\mathcal D^\perp)$ is the kernel of the pulled-back K\"ahler form. Here $F_*^{-1}(\mathcal D^\perp)$ denotes the full preimage in $TM$, including $\mathcal V$. Without the maximality hypothesis, the chosen summand need not be integrable: for instance, a totally real submanifold can admit a nonintegrable proper subdistribution chosen as $\mathcal D^\perp$.

\begin{proposition}\label{prop:ptr_characteristic}
Let $F:(M,g_M)\to(N,g_N,J)$ be a PTR Riemannian map satisfying $\ker P=\mathcal D^\perp$. Then $\dim\mathcal D$ is even and $F_*^{-1}(\mathcal D^\perp)$ is integrable. Moreover, its horizontal part is integrable if and only if
\begin{equation}\label{eq:ptr_characteristic_A}
\mathcal A_WZ=0
\end{equation}
for all horizontal vector fields $W,Z$ such that $F_*W,F_*Z\in\mathcal D^\perp$.
\end{proposition}

\begin{proof}
Let $\Omega$ be the K\"ahler form on $N$ and put $\omega=F^*\Omega$. Since $d\Omega=0$, we have $d\omega=0$. For $X,Y\in\Gamma(TM)$,
\[
\omega(X,Y)=g_N(JF_*X,F_*Y).
\]
Consequently, $\ker\omega=F_*^{-1}(\ker P)=F_*^{-1}(\mathcal D^\perp)$. Its rank is constant by the PTR hypothesis, and the kernel of a closed $2$-form of constant rank is integrable. Since $P$ is skew-adjoint and has trivial kernel on $\mathcal D$, the dimension of $\mathcal D$ is even.

If $W,Z$ are horizontal and $F_*W,F_*Z\in\mathcal D^\perp$, the first assertion yields $[W,Z]\in\Gamma(F_*^{-1}(\mathcal D^\perp))$. Thus the horizontal component of $[W,Z]$ lies in the horizontal lift of $\mathcal D^\perp$. Since $\mathcal V[W,Z]=2\mathcal A_WZ$, this horizontal lift is integrable precisely when \eqref{eq:ptr_characteristic_A} holds.
\end{proof}

We next consider the horizontal lift of $\mathcal D$. Its integrability and total geodesicity are characterized without assuming maximality of $\mathcal D^\perp$; compare the corresponding distributional conditions for PTR submanifolds~\cite{PoyrazSahinYerlikaya2025}.

\begin{proposition}\label{prop:ptr_ambiguous_geometry}
Let $F:(M,g_M)\to(N,g_N,J)$ be a PTR Riemannian map. For horizontal vector fields $X,Y$ satisfying $F_*X,F_*Y\in\mathcal D$, let $X_P,Y_P$ be the unique horizontal vector fields determined by $F_*X_P=PF_*X$ and $F_*Y_P=PF_*Y$. Then the following statements hold.
\begin{enumerate}[label=(\roman*)]
\item The horizontal lift of $\mathcal D$ is integrable if and only if $\mathcal A_XY=0$ and
\begin{equation}\label{eq:ptr_ambiguous_integrability}
\begin{aligned}
0=g_N\big(& (\nabla F_*)(X,Y_P)-(\nabla F_*)(Y,X_P) \\
&+\nabla_X^{F,\perp}(QF_*Y)-\nabla_Y^{F,\perp}(QF_*X),JF_*W\big)
\end{aligned}
\end{equation}
for every horizontal $W$ such that $F_*W\in\mathcal D^\perp$.
\item The horizontal lift of $\mathcal D$ defines a totally geodesic foliation on $M$ if and only if $\mathcal A_XY=0$ and
\begin{equation}\label{eq:ptr_ambiguous_geodesic}
g_N((\nabla F_*)(X,Y_P)+\nabla_X^{F,\perp}(QF_*Y),JF_*W)=0
\end{equation}
for all $X,Y,W$ as above.
\end{enumerate}
\end{proposition}

\begin{proof}
For the horizontal vector fields in the statement, the K\"ahler identity and the Gauss--Weingarten formulas give
\begin{equation}\label{eq:ptr_ambiguous_connection}
g_M(\nabla_X^M Y,W)=g_N((\nabla F_*)(X,Y_P)+\nabla_X^{F,\perp}(QF_*Y),JF_*W).
\end{equation}
To see this, the normality of $(\nabla F_*)(X,Y)$ implies
\[
g_M(\nabla_X^M Y,W)=g_N(\nabla_X^F F_*Y,F_*W)=g_N(\nabla_X^F(JF_*Y),JF_*W).
\]
Now use $JF_*Y=F_*Y_P+QF_*Y$ and take the normal component of the covariant derivative on the right. This proves \eqref{eq:ptr_ambiguous_connection}.

Interchanging $X$ and $Y$ in \eqref{eq:ptr_ambiguous_connection} and subtracting gives the component of $[X,Y]$ along the horizontal lift of $\mathcal D^\perp$. The remaining obstruction to integrability is its vertical component $\mathcal V[X,Y]=2\mathcal A_XY$. This proves (i).

For (ii), the horizontal lift of $\mathcal D$ defines a totally geodesic foliation precisely when $\nabla_X^M Y$ belongs to that distribution for all such $X,Y$. Its vertical component is $\mathcal A_XY$, and its component along the horizontal lift of $\mathcal D^\perp$ is given by \eqref{eq:ptr_ambiguous_connection}. Their simultaneous vanishing is equivalent to the stated conditions.
\end{proof}

The next proposition gives a totally geodesic criterion for the horizontal lift of $\mathcal D^\perp$. Its stated form uses $\ker P=\mathcal D^\perp$: without this condition, the vanishing of $PF_*(\nabla^M_WZ)$ need not eliminate the component along $\mathcal D$. More general criteria can be formulated by retaining that component.

\begin{proposition}\label{prop:ptr_totally_real_geodesic}
Let $F:(M,g_M)\to(N,g_N,J)$ be a PTR Riemannian map satisfying $\ker P=\mathcal D^\perp$. The horizontal lift of $\mathcal D^\perp$ defines a totally geodesic foliation on $M$ if and only if
\begin{equation}\label{eq:ptr_totally_real_geodesic}
\mathcal A_WZ=0,\qquad \mathcal S_{JF_*Z}F_*W+t(\nabla F_*)(W,Z)=0
\end{equation}
for all horizontal vector fields $W,Z$ with $F_*W,F_*Z\in\mathcal D^\perp$.
\end{proposition}

\begin{proof}
The horizontal lift of $\mathcal D^\perp$ defines a totally geodesic foliation if and only if $\nabla_W^M Z$ has no component in $\mathcal V$ and no component in the horizontal lift of $\mathcal D$. The first condition is $\mathcal A_WZ=0$.

To examine the second, use $JF_*Z\in(\operatorname{range}F_*)^\perp$. The K\"ahler identity and the Gauss--Weingarten formulas yield
\[
-\mathcal S_{JF_*Z}F_*W+\nabla_W^{F,\perp}(JF_*Z)=JF_*(\nabla_W^MZ)+J(\nabla F_*)(W,Z).
\]
Comparing the components tangent to $\operatorname{range}F_*$, we find
\[
PF_*(\nabla_W^MZ)=-\mathcal S_{JF_*Z}F_*W-t(\nabla F_*)(W,Z).
\]
Since $\ker P=\mathcal D^\perp$, the right-hand side vanishes exactly when the horizontal component of $\nabla_W^MZ$ projects into $\mathcal D^\perp$. Together with $\mathcal A_WZ=0$, this proves both implications.
\end{proof}

When the horizontal distribution is integrable, the same maximality condition leads to a local product criterion in terms of the shape operator and $(\nabla F_*)$, as in the corresponding PTR-submanifold setting~\cite{PoyrazSahinYerlikaya2025}.

\begin{theorem}\label{thm:ptr_horizontal_product}
Let $F:(M,g_M)\to(N,g_N,J)$ be a PTR Riemannian map satisfying $\ker P=\mathcal D^\perp$, and suppose that $\mathcal H$ is integrable. Then every leaf of $\mathcal H$ is locally a Riemannian product of the integral distributions corresponding to $\mathcal D$ and $\mathcal D^\perp$ if and only if
\begin{equation}\label{eq:ptr_horizontal_product}
\mathcal S_{JF_*W}F_*X+t(\nabla F_*)(X,W)=0
\end{equation}
for all horizontal vector fields $X,W$ such that $F_*W\in\mathcal D^\perp$.
\end{theorem}

\begin{proof}
Since $JF_*W$ is normal, the Weingarten formula and the K\"ahler condition imply
\[
-\mathcal S_{JF_*W}F_*X+\nabla_X^{F,\perp}(JF_*W)=JF_*(\nabla_X^MW)+J(\nabla F_*)(X,W).
\]
Taking the component tangent to $\operatorname{range}F_*$ gives
\begin{equation}\label{eq:ptr_product_identity}
PF_*(\nabla_X^MW)=-\mathcal S_{JF_*W}F_*X-t(\nabla F_*)(X,W).
\end{equation}
Because $\ker P=\mathcal D^\perp$, equation \eqref{eq:ptr_horizontal_product} is equivalent to the requirement that the horizontal component of $\nabla_X^MW$ project into $\mathcal D^\perp$ for every horizontal $X$. On an integral leaf of $\mathcal H$, this means that the distribution corresponding to $\mathcal D^\perp$ is parallel for the induced Levi-Civita connection. Its orthogonal complement, corresponding to $\mathcal D$, is therefore parallel as well. The local Riemannian product decomposition follows. Conversely, if the two factors are parallel, \eqref{eq:ptr_product_identity} gives \eqref{eq:ptr_horizontal_product}.
\end{proof}

Finally, we consider horizontally totally umbilical PTR Riemannian maps. Let $r=\operatorname{rank}F$, and denote by
\[
H_F=\frac1r\sum_{a=1}^r(\nabla F_*)(E_a,E_a)
\]
the mean curvature vector of the range distribution, where $\{E_1,\ldots,E_r\}$ is any orthonormal frame of $\mathcal H$. We obtain a restriction on $H_F$ which will be relevant to the inequality in Section~4.

\begin{proposition}\label{prop:ptr_umbilical_mean}
Let $F:(M,g_M)\to(N,g_N,J)$ be a PTR Riemannian map with $q=\dim\mathcal D^\perp\geq2$. Suppose that
\begin{equation}\label{eq:ptr_horizontal_umbilical}
(\nabla F_*)(X,Y)=g_M(X,Y)H_F,\qquad X,Y\in\Gamma(\mathcal H).
\end{equation}
Then $H_F$ is orthogonal to $J\mathcal D^\perp$, that is,
\begin{equation}\label{eq:ptr_umbilical_orthogonal}
g_N(H_F,JF_*Z)=0
\end{equation}
for every horizontal $Z$ with $F_*Z\in\mathcal D^\perp$.
\end{proposition}

\begin{proof}
Let $Z$ be a unit horizontal vector with $F_*Z\in\mathcal D^\perp$. Because $q\geq2$, we can choose a unit horizontal vector $W$ with $F_*W\in\mathcal D^\perp$ and $g_M(W,Z)=0$. Lemma~\ref{lem:ptr_cubic} and \eqref{eq:ptr_horizontal_umbilical} give
\[
\begin{aligned}
g_N(H_F,JF_*Z)&=g_N((\nabla F_*)(W,W),JF_*Z)\\
&=g_N((\nabla F_*)(W,Z),JF_*W)=0.
\end{aligned}
\]
The assertion follows by linearity.
\end{proof}

\begin{corollary}\label{cor:ptr_umbilical_tot_geodesic}
Under the hypotheses of Proposition~\ref{prop:ptr_umbilical_mean}, suppose also that $(\operatorname{range}F_*)^\perp=J\mathcal D^\perp$. Then the horizontal second fundamental form vanishes.
\end{corollary}

\begin{proof}
Since $H_F$ is normal to $\operatorname{range}F_*$ and orthogonal to $J\mathcal D^\perp$, the additional hypothesis implies $H_F=0$. The conclusion follows from \eqref{eq:ptr_horizontal_umbilical}.
\end{proof}

There is a further consequence for harmonic maps, which also clarifies the distinction between the vanishing of the horizontal second fundamental form and total geodesicity of the map. Recall that the tension field of a Riemannian map contains independent horizontal and vertical contributions.

\begin{corollary}\label{cor:ptr_umbilical_harmonic}
Let $F:(M,g_M)\to(N,g_N,J)$ be a horizontally totally umbilical PTR Riemannian map satisfying $q\geq2$ and $(\operatorname{range}F_*)^\perp=J\mathcal D^\perp$. Then $F$ is harmonic if and only if its fibers are minimal.
\end{corollary}

\begin{proof}
By Corollary~\ref{cor:ptr_umbilical_tot_geodesic}, the horizontal second fundamental form vanishes and $H_F=0$. Let $V_1,\ldots,V_{m-r}$ be an orthonormal frame of $\mathcal V$, where $m=\dim M$. The tension field decomposes as
\[
\tau(F)=-F_*\left(\sum_{a=1}^{m-r}\mathcal H\nabla_{V_a}^M V_a\right)+rH_F.
\]
Indeed, $(\nabla F_*)(V_a,V_a)=-F_*(\nabla_{V_a}^MV_a)$, while the trace of $(\nabla F_*)$ over $\mathcal H$ is $rH_F$. Hence, under the hypotheses, $\tau(F)=0$ is equivalent to the vanishing of the mean curvature of the fibers.
\end{proof}

\begin{remark}
The preceding mean-curvature restriction is related to the results on totally umbilical PTR submanifolds in \cite{PoyrazSahinYerlikaya2025}, while the corresponding notion for Riemannian maps is discussed in \cite{Sahin2011}. Its formulation here uses the horizontal second fundamental form of a Riemannian map. In particular, for $q\geq2$ the component of $H_F$ along $J\mathcal D^\perp$ vanishes in the totally umbilical case, so the corresponding correction term in Theorem~\ref{thm:ptr_refined_delta} is zero.
\end{remark}

\section{Curvature Inequalities for PTR Riemannian Maps}

Let $F:(M,g_M)\to (N(c),g_N,J)$ be a PTR Riemannian map into a complex space form of constant holomorphic sectional curvature $c$. Suppose that
\[
\dim\mathcal D=d,\qquad \dim\mathcal D^\perp=q,\qquad \operatorname{rank}F=r=d+q.
\]
At $p\in M$, choose an orthonormal basis
\[
\{E_1,\ldots,E_r\}=\{X_1,\ldots,X_d,W_1,\ldots,W_q\}
\]
of $\mathcal H_p$ such that $F_*X_i\in\mathcal D$ and $F_*W_\alpha\in\mathcal D^\perp$.

We denote the horizontal scalar curvature of $M$ by
\[
\tau^{\mathcal H}=\sum_{1\leq a<b\leq r}K_M(E_a,E_b),
\]
and define
\[
\tau^{\mathcal H}(\mathcal D)=\sum_{1\leq i<j\leq d}K_M(X_i,X_j).
\]
These are traces of the sectional curvature of $M$ over the corresponding horizontal subspaces. In analogy with the CR $\delta$-invariant~\cite{Chen2012CR} and the PTR curvature trace~\cite{MihaiEtAl2026}, we define
\begin{equation}\label{eq:ptr_delta}
\delta_{\mathrm{PTR}}^{\mathcal H}(\mathcal D) =\tau^{\mathcal H}-\tau^{\mathcal H}(\mathcal D).
\end{equation}
The quantity in \eqref{eq:ptr_delta} is determined by the chosen PTR decomposition and need not be canonical when $\mathcal D^\perp$ is not maximal. In the adapted frame,
\begin{equation}\label{eq:ptr_delta_expansion}
\delta_{\mathrm{PTR}}^{\mathcal H}(\mathcal D) =\sum_{i=1}^d\sum_{\alpha=1}^q K_M(X_i,W_\alpha) +\sum_{1\leq\alpha<\beta\leq q}K_M(W_\alpha,W_\beta).
\end{equation}

In this adapted frame, the mean curvature vector of the range distribution has the form
\[
H_F=\frac1r\left\{ \sum_{i=1}^d(\nabla F_*)(X_i,X_i) +\sum_{\alpha=1}^q(\nabla F_*)(W_\alpha,W_\alpha) \right\}.
\]

We first establish an identity which will be used to derive the curvature inequality.

\begin{proposition}\label{prop:ptr_curvature_identity}
Let $F:(M,g_M)\to(N(c),g_N,J)$ be a PTR Riemannian map of rank $r=d+q$. Then
\begin{equation}\label{eq:ptr_curvature_identity}
\begin{aligned}
\delta_{\mathrm{PTR}}^{\mathcal H}(\mathcal D)
={}&\frac{c}{8}q(2d+q-1)+\frac{r^2}{2}\|H_F\|^2\\
&-\frac12\left\|\sum_{i=1}^d(\nabla F_*)(X_i,X_i)\right\|^2
-\frac12\sum_{\alpha,\beta=1}^q\|(\nabla F_*)(W_\alpha,W_\beta)\|^2\\
&-\sum_{i=1}^d\sum_{\alpha=1}^q\|(\nabla F_*)(X_i,W_\alpha)\|^2\\
&-3\sum_{i=1}^d\sum_{\alpha=1}^q\|\mathcal A_{X_i}W_\alpha\|^2
-3\sum_{1\leq\alpha<\beta\leq q}\|\mathcal A_{W_\alpha}W_\beta\|^2.
\end{aligned}
\end{equation}
\end{proposition}

\begin{proof}
For orthonormal horizontal vectors $X,Y$, the Gauss equation for Riemannian maps, together with O'Neill's horizontal curvature formula~\cite{ONeill1966,Sahin2017Book}, gives
\begin{equation}\label{eq:ptr_gauss_oneill}
\begin{aligned}
K_M(X,Y)={}&K_N(F_*X,F_*Y)
+g_N((\nabla F_*)(X,X),(\nabla F_*)(Y,Y))\\
&-\|(\nabla F_*)(X,Y)\|^2-3\|\mathcal A_XY\|^2.
\end{aligned}
\end{equation}
Since $\mathcal D^\perp$ is totally real, we have
\[
g_N(JF_*X_i,F_*W_\alpha)=0,\qquad g_N(JF_*W_\alpha,F_*W_\beta)=0.
\]
Consequently, the curvature tensor of $N(c)$ yields
\[
K_N(F_*X_i,F_*W_\alpha) =K_N(F_*W_\alpha,F_*W_\beta)=\frac c4
\]
for distinct $\alpha,\beta$.

Using \eqref{eq:ptr_delta_expansion} and \eqref{eq:ptr_gauss_oneill}, we obtain
\begin{align*}
\delta_{\mathrm{PTR}}^{\mathcal H}(\mathcal D)
={}&\frac c8q(2d+q-1)\\
&+\sum_{i,\alpha}g_N((\nabla F_*)(X_i,X_i),(\nabla F_*)(W_\alpha,W_\alpha))\\
&+\sum_{\alpha<\beta}g_N((\nabla F_*)(W_\alpha,W_\alpha),(\nabla F_*)(W_\beta,W_\beta))\\
&-\sum_{i,\alpha}\|(\nabla F_*)(X_i,W_\alpha)\|^2
-\sum_{\alpha<\beta}\|(\nabla F_*)(W_\alpha,W_\beta)\|^2\\
&-3\sum_{i,\alpha}\|\mathcal A_{X_i}W_\alpha\|^2
-3\sum_{\alpha<\beta}\|\mathcal A_{W_\alpha}W_\beta\|^2.
\end{align*}

By the definition of $H_F$,
\begin{align*}
&\sum_{i,\alpha}g_N((\nabla F_*)(X_i,X_i),(\nabla F_*)(W_\alpha,W_\alpha))\\
&\quad+\sum_{\alpha<\beta}g_N((\nabla F_*)(W_\alpha,W_\alpha),(\nabla F_*)(W_\beta,W_\beta))\\
&=\frac{r^2}{2}\|H_F\|^2
-\frac12\left\|\sum_i(\nabla F_*)(X_i,X_i)\right\|^2
-\frac12\sum_\alpha\|(\nabla F_*)(W_\alpha,W_\alpha)\|^2.
\end{align*}
Combining this identity with
\[
\sum_{\alpha<\beta}\|(\nabla F_*)(W_\alpha,W_\beta)\|^2 +\frac12\sum_\alpha\|(\nabla F_*)(W_\alpha,W_\alpha)\|^2 =\frac12\sum_{\alpha,\beta}\|(\nabla F_*)(W_\alpha,W_\beta)\|^2
\]
gives \eqref{eq:ptr_curvature_identity}.
\end{proof}

We now obtain a refined estimate by using the symmetry of the normal components associated with the totally real distribution.

\begin{theorem}\label{thm:ptr_refined_delta}
Let $F:(M,g_M)\to(N(c),g_N,J)$ be a PTR Riemannian map of rank $r=d+q$, where $d=\dim\mathcal D\geq1$ and $q=\dim\mathcal D^\perp\geq1$. Then
\begin{equation}\label{eq:ptr_refined_delta}
\begin{aligned}
\delta_{\mathrm{PTR}}^{\mathcal H}(\mathcal D)
\leq{}&\frac{c}{8}q(2d+q-1)
+\frac{r^2q}{2(q+1)}\|H_F\|^2\\
&-\frac{r^2(q-1)}{(q+1)(q+5)}
\sum_{\alpha=1}^q g_N(H_F,JF_*W_\alpha)^2\\
&-\sum_{i=1}^d\sum_{\alpha=1}^q\|(\nabla F_*)(X_i,W_\alpha)\|^2\\
&-3\sum_{i=1}^d\sum_{\alpha=1}^q\|\mathcal A_{X_i}W_\alpha\|^2
-3\sum_{1\leq\alpha<\beta\leq q}\|\mathcal A_{W_\alpha}W_\beta\|^2.
\end{aligned}
\end{equation}
\end{theorem}

\begin{proof}
Fix $p\in M$ and retain the adapted orthonormal basis $\{E_1,\ldots,E_r\}$. Choose an orthonormal basis $\{\xi_1,\ldots,\xi_\ell\}$ of $(\operatorname{range}F_*)^\perp$ such that
\[
\xi_\alpha=JF_*W_\alpha,\qquad 1\leq\alpha\leq q.
\]
Write
\[
h_{ab}^{\lambda}=g_N((\nabla F_*)(E_a,E_b),\xi_\lambda), \qquad 1\leq a,b\leq r,\quad 1\leq\lambda\leq\ell,
\]
and put
\[
u_\lambda=\sum_{\alpha=1}^q h_{d+\alpha,d+\alpha}^{\lambda}, \qquad v_\lambda=\sum_{i=1}^d h_{ii}^{\lambda}.
\]
Then
\begin{equation}\label{eq:ptr_mean_components}
u_\lambda+v_\lambda=r\,g_N(H_F,\xi_\lambda).
\end{equation}

By Proposition~\ref{prop:ptr_curvature_identity},
\begin{equation}\label{eq:ptr_refined_step1}
\begin{aligned}
\delta_{\mathrm{PTR}}^{\mathcal H}(\mathcal D)
={}&\frac{c}{8}q(2d+q-1)+\frac{r^2}{2}\|H_F\|^2\\
&-\frac12\sum_{\lambda=1}^{\ell}
\left\{v_\lambda^2+\sum_{\alpha,\beta=1}^q
(h_{d+\alpha,d+\beta}^{\lambda})^2\right\}\\
&-\sum_{i,\alpha}\|(\nabla F_*)(X_i,W_\alpha)\|^2
-3\sum_{i,\alpha}\|\mathcal A_{X_i}W_\alpha\|^2\\
&-3\sum_{\alpha<\beta}\|\mathcal A_{W_\alpha}W_\beta\|^2.
\end{aligned}
\end{equation}

By Lemma~\ref{lem:ptr_cubic}, the components
\[
h_{d+\alpha,d+\beta}^{\gamma}, \qquad 1\leq\alpha,\beta,\gamma\leq q,
\]
are symmetric in all three indices. The associated cubic estimate, used also in the PTR-submanifold setting~\cite{MihaiEtAl2026}, gives
\begin{equation}\label{eq:ptr_cubic_estimate}
\sum_{\alpha,\beta,\gamma=1}^q(h_{d+\alpha,d+\beta}^{\gamma})^2 \geq\frac{3}{q+2}\sum_{\gamma=1}^q u_\gamma^2.
\end{equation}
Indeed, the difference between the two sides of \eqref{eq:ptr_cubic_estimate} equals
\[
\sum_{\alpha,\beta,\gamma=1}^q \left( h_{d+\alpha,d+\beta}^{\gamma} -\frac{\delta_{\alpha\beta}u_\gamma+ \delta_{\alpha\gamma}u_\beta+ \delta_{\beta\gamma}u_\alpha}{q+2} \right)^2,
\]
which is nonnegative.

For the remaining normal directions, the Cauchy--Schwarz inequality gives
\begin{equation}\label{eq:ptr_remaining_estimate}
\sum_{\alpha,\beta=1}^q(h_{d+\alpha,d+\beta}^{\lambda})^2 \geq\frac1q u_\lambda^2,\qquad q<\lambda\leq\ell.
\end{equation}

Furthermore,
\begin{equation}\label{eq:ptr_trace_estimates}
\begin{aligned}
v_\gamma^2+\frac{3}{q+2}u_\gamma^2
&\geq\frac{3}{q+5}(u_\gamma+v_\gamma)^2,
&&1\leq\gamma\leq q,\\
v_\lambda^2+\frac1q u_\lambda^2
&\geq\frac1{q+1}(u_\lambda+v_\lambda)^2,
&&q<\lambda\leq\ell.
\end{aligned}
\end{equation}

Applying \eqref{eq:ptr_cubic_estimate}--\eqref{eq:ptr_trace_estimates} to \eqref{eq:ptr_refined_step1}, and using \eqref{eq:ptr_mean_components}, we obtain
\begin{align*}
\delta_{\mathrm{PTR}}^{\mathcal H}(\mathcal D)
\leq{}&\frac{c}{8}q(2d+q-1)
+\frac{r^2(q+2)}{2(q+5)}
\sum_{\gamma=1}^q g_N(H_F,\xi_\gamma)^2\\
&+\frac{r^2q}{2(q+1)}
\sum_{\lambda=q+1}^{\ell}g_N(H_F,\xi_\lambda)^2\\
&-\sum_{i,\alpha}\|(\nabla F_*)(X_i,W_\alpha)\|^2
-3\sum_{i,\alpha}\|\mathcal A_{X_i}W_\alpha\|^2\\
&-3\sum_{\alpha<\beta}\|\mathcal A_{W_\alpha}W_\beta\|^2.
\end{align*}

Finally, since
\[
\sum_{\lambda=1}^{\ell}g_N(H_F,\xi_\lambda)^2=\|H_F\|^2
\]
and
\[
\frac{q+2}{q+5}-\frac{q}{q+1} =-\frac{2(q-1)}{(q+1)(q+5)},
\]
the preceding inequality reduces to \eqref{eq:ptr_refined_delta}.
\end{proof}

\noindent\textit{Equality case.}
Equality in \eqref{eq:ptr_refined_delta} at $p\in M$ holds if and only if, with respect to an adapted orthonormal frame as in the proof,
\begin{equation}\label{eq:ptr_equality_one}
h_{d+\alpha,d+\beta}^{\gamma} =\frac{\delta_{\alpha\beta}u_\gamma+ \delta_{\alpha\gamma}u_\beta+ \delta_{\beta\gamma}u_\alpha}{q+2}, \qquad 1\leq\alpha,\beta,\gamma\leq q,
\end{equation}
\begin{equation}\label{eq:ptr_equality_two}
h_{d+\alpha,d+\beta}^{\lambda} =\frac{u_\lambda}{q}\delta_{\alpha\beta}, \qquad q<\lambda\leq\ell,
\end{equation}
and
\begin{equation}\label{eq:ptr_equality_three}
v_\gamma=\frac{3}{q+2}u_\gamma,\qquad v_\lambda=\frac1q u_\lambda,
\end{equation}
for $1\leq\gamma\leq q$ and $q<\lambda\leq\ell$, respectively.

Indeed, \eqref{eq:ptr_equality_one} and \eqref{eq:ptr_equality_two} characterize equality in \eqref{eq:ptr_cubic_estimate} and \eqref{eq:ptr_remaining_estimate}, while \eqref{eq:ptr_equality_three} is equivalent to equality in \eqref{eq:ptr_trace_estimates}. Since the remaining terms in \eqref{eq:ptr_refined_delta} are retained, equality imposes no additional vanishing condition on the mixed second fundamental form or the tensor $\mathcal A$.

\begin{proposition}\label{prop:ptr_mihai}
Let $F:(M,g_M)\to(N(c),g_N,J)$ be a PTR Riemannian map with $d=\dim\mathcal D\geq1$, $q=\dim\mathcal D^\perp\geq1$ and $r=d+q$. Put
\[
\nu=(\operatorname{range}F_*)^\perp\ominus J\mathcal D^\perp
\]
and denote by $H_{\mathcal D^\perp}^{J}$ and $H_{\mathcal D^\perp}^{\nu}$ the orthogonal projections of $H_{\mathcal D^\perp}$ onto $J\mathcal D^\perp$ and $\nu$, respectively. Then
\begin{equation}\label{eq:ptr_mihai}
\begin{aligned}
\delta_{\mathrm{PTR}}^{\mathcal H}(\mathcal D)\leq{}&
\frac{c}{8}q(2d+q-1)+\frac{r^2}{2}\|H_F\|^2
-\frac{d^2}{2}\|H_{\mathcal D}\|^2\\
&-\frac{3q^2}{2(q+2)}\|H_{\mathcal D^\perp}^{J}\|^2
-\frac q2\|H_{\mathcal D^\perp}^{\nu}\|^2\\
&-\sum_{i=1}^d\sum_{\alpha=1}^q\|(\nabla F_*)(X_i,W_\alpha)\|^2
-3\sum_{i=1}^d\sum_{\alpha=1}^q\|\mathcal A_{X_i}W_\alpha\|^2\\
&-3\sum_{1\leq\alpha<\beta\leq q}\|\mathcal A_{W_\alpha}W_\beta\|^2,
\end{aligned}
\end{equation}
where
\[
H_{\mathcal D}=\frac1d\sum_{i=1}^d(\nabla F_*)(X_i,X_i),\qquad H_{\mathcal D^\perp}=\frac1q\sum_{\alpha=1}^q(\nabla F_*)(W_\alpha,W_\alpha).
\]
Equality at a point holds if and only if \eqref{eq:ptr_equality_one} and \eqref{eq:ptr_equality_two} hold with respect to an adapted orthonormal frame.
\end{proposition}

\begin{proof}
Complete $\{JF_*W_1,\ldots,JF_*W_q\}$ to an orthonormal frame of $(\operatorname{range}F_*)^\perp$. The cubic estimate \eqref{eq:ptr_cubic_estimate} on $J\mathcal D^\perp$ and the Cauchy--Schwarz estimate \eqref{eq:ptr_remaining_estimate} on $\nu$ give
\[
\sum_{\alpha,\beta=1}^q\|(\nabla F_*)(W_\alpha,W_\beta)\|^2\geq\frac{3q^2}{q+2}\|H_{\mathcal D^\perp}^{J}\|^2+q\|H_{\mathcal D^\perp}^{\nu}\|^2.
\]
The first partial trace in \eqref{eq:ptr_curvature_identity} is $dH_{\mathcal D}$. Substitution into that identity gives \eqref{eq:ptr_mihai}. Equality holds if and only if equality occurs in both \eqref{eq:ptr_cubic_estimate} and \eqref{eq:ptr_remaining_estimate}, namely, if and only if \eqref{eq:ptr_equality_one} and \eqref{eq:ptr_equality_two} hold.
\end{proof}

\begin{remark}
The estimate \eqref{eq:ptr_mihai} retains the partial mean curvatures and is valid for proper PTR Riemannian maps in arbitrary normal codimension. Optimizing these partial traces, as in \eqref{eq:ptr_trace_estimates}, gives Theorem~\ref{thm:ptr_refined_delta}, which is expressed in terms of $H_F$ and the projection of $H_F$ onto $J\mathcal D^\perp$.

If $(\operatorname{range}F_*)^\perp=J\mathcal D^\perp$, then $\nu=\{0\}$. Since $Q\mathcal D$ is orthogonal to $J\mathcal D^\perp$ by the calculation in Lemma~\ref{lem:ptr_normal_decomposition}, we have $Q\mathcal D=0$; thus $J\mathcal D=\mathcal D$ and $d$ is even. In this minimum-codimension case, omitting the nonpositive $H_{\mathcal D}$, mixed, and O'Neill terms recovers the inequality of Mihai et al.~\cite[Theorem~1]{MihaiEtAl2026} for an isometric immersion. The same condition makes the corresponding PTR submanifold a $CR$-submanifold.
\end{remark}

The following example shows that equality can occur for a proper PTR Riemannian map with nonzero mean curvature and nonintegrable horizontal distribution.

\begin{example}\label{ex:ptr_sharp}
Let $d\geq4$, $q\geq1$, and $r=d+q$. Choose $0<t_1<t_2<1$ and put $\Lambda_j=\sqrt{1+t_j^2}$, $j=1,2$. Define $L:\mathbb R^r\to\mathbb C^r$ by
\[
\begin{aligned}
L(s,w)=\bigg(&\frac{s_1-it_1s_2}{\Lambda_1},
\frac{s_2+it_1s_1}{\Lambda_1},
\frac{s_3-it_2s_4}{\Lambda_2},\\
&\frac{s_4+it_2s_3}{\Lambda_2},
s_5,\ldots,s_d,w_1,\ldots,w_q\bigg),
\end{aligned}
\]
where $s=(s_1,\ldots,s_d)$ and $w=(w_1,\ldots,w_q)$.

For $x=(x_1,\ldots,x_r)\in\mathbb R^r$, set
\[
S(x)=\Lambda_1^2(x_1^2+x_2^2)+\Lambda_2^2(x_3^2+x_4^2)+\sum_{j=5}^d x_j^2
\]
and
\[
\begin{aligned}
K(x)={}&\frac{\varepsilon}{2(q+2)}x_{d+1}\sum_{\alpha=1}^q x_{d+\alpha}^2
+\frac{3\varepsilon}{2d(q+2)}x_{d+1}S(x)\\
&+\frac{\kappa}{2}x_1\sum_{\alpha=1}^q x_{d+\alpha}^2
+\frac{\kappa}{2(d+2)}x_1S(x),
\end{aligned}
\]
where $\varepsilon\kappa\ne0$.

Writing $z=x+iy\in\mathbb C^r$, define
\[
\Phi(x+iy)=x+i\bigl(y+\nabla K(x)\bigr),\qquad G=\Phi\circ L.
\]
Let $U$ be a sufficiently small neighborhood of the origin in $\mathbb R^r$. On $M=U\times\mathbb R$, define
\[
F(s,w,v)=G(s,w)
\]
and
\[
g_M=F^*g_{\mathbb C^r}+(dv+s_1\,dw_1)^2.
\]
Since $G$ is an immersion near the origin, $g_M$ is a Riemannian metric and $F$ is a Riemannian map with $\ker F_*=\operatorname{span}\{\partial_v\}$.

The transformation $\Phi$ preserves the standard K\"ahler form $\Omega$, since the Hessian of $K$ is symmetric. Hence
\[
G^*\Omega=L^*\Omega =-\frac{2t_1}{1+t_1^2}\,ds_1\wedge ds_2 -\frac{2t_2}{1+t_2^2}\,ds_3\wedge ds_4.
\]
Consequently, the images of the horizontal lifts of $\partial_{w_1},\ldots,\partial_{w_q}$ span a totally real distribution $\mathcal D^\perp$. Its orthogonal complement $\mathcal D$ in $\operatorname{range}F_*$ is preserved by $P$.

At the origin, $P$ has two distinct nonzero K\"ahler-angle parameters,
\[
\frac{2t_1}{1+t_1^2},\qquad \frac{2t_2}{1+t_2^2}.
\]
Thus $\mathcal D$ is neither invariant nor pointwise slant. By continuity, $F$ is a proper PTR Riemannian map on a sufficiently small neighborhood. When $d=4$, the chosen totally real summand is maximal. For $d\geq5$, the directions corresponding to $s_5,\ldots,s_d$ also belong to $\ker P$, so $\ker P\supsetneq\mathcal D^\perp$; this is permitted by Definition~\ref{def:ptr_map}.

We verify equality at the origin. There the horizontal vectors
\[
X_i=\partial_{s_i}\quad(1\leq i\leq d),\qquad W_1=\partial_{w_1}-s_1\partial_v,\qquad W_\alpha=\partial_{w_\alpha}\quad(2\leq\alpha\leq q)
\]
form an orthonormal basis. Let
\[
\zeta=\frac{t_1\partial_{x_2}+\partial_{y_1}}{1+t_1^2}.
\]
Then $\zeta$ is a nonzero normal vector orthogonal to $J\mathcal D^\perp$.

Direct calculation of the normal components of the second derivatives of $F$ gives
\[
\begin{aligned}
(\nabla F_*)(W_\alpha,W_\beta)
={}&\frac{\varepsilon}{q+2}\bigl(
\delta_{\alpha\beta}JF_*W_1
+\delta_{\alpha1}JF_*W_\beta
+\delta_{\beta1}JF_*W_\alpha\bigr)\\
&+\kappa\delta_{\alpha\beta}\zeta,
\end{aligned}
\]
and
\[
\sum_{i=1}^d(\nabla F_*)(X_i,X_i) =\frac{3\varepsilon}{q+2}JF_*W_1+\kappa\zeta.
\]
It follows that
\[
\sum_{\alpha=1}^q(\nabla F_*)(W_\alpha,W_\alpha) =\varepsilon JF_*W_1+q\kappa\zeta.
\]
Therefore,
\[
H_F=\frac{1}{r}\left( \frac{q+5}{q+2}\varepsilon JF_*W_1 +(q+1)\kappa\zeta \right).
\]

The preceding expressions satisfy the cubic symmetry equality condition \eqref{eq:ptr_equality_one}. In the normal directions orthogonal to $J\mathcal D^\perp$, the restriction of $(\nabla F_*)$ to the span of $W_1,\ldots,W_q$ is a scalar multiple of the metric. Moreover, the two partial traces satisfy precisely the relations in \eqref{eq:ptr_equality_three}. Thus equality holds in \eqref{eq:ptr_refined_delta} at the origin. In particular, the conditions \eqref{eq:ptr_equality_one} and \eqref{eq:ptr_equality_two} are satisfied, so equality is also attained in Proposition~\ref{prop:ptr_mihai}.

Finally,
\[
[X_1,W_1]=-\partial_v,\qquad \mathcal A_{X_1}W_1=-\frac12\partial_v\ne0.
\]
Hence the horizontal distribution is nonintegrable. Since $\varepsilon\kappa\ne0$, both the component of $H_F$ along $J\mathcal D^\perp$ and its orthogonal normal component are nonzero.

For $q\geq2$, equality occurs with nonzero components of $H_F$ both along $J\mathcal D^\perp$ and orthogonal to it. Thus neither coefficient in the two-component mean-curvature estimate in the proof of Theorem~\ref{thm:ptr_refined_delta} can be decreased independently while the other coefficients remain unchanged. For $q=1$, the correction term in \eqref{eq:ptr_refined_delta} vanishes and the remaining mean-curvature coefficient is attained.
\end{example}

\end{document}